\documentclass{owrart}

\usepackage{amssymb}

\newcommand{\PP}{{\mathbb P}}
\newcommand{\Q}{{\mathbb Q}}
\newcommand{\Qdot}{{\dot\Q}}
\newcommand{\of}{\subseteq}
\newcommand{\ofnoteq}{\subsetneq}
\newcommand{\ZFC}{\mathsf{ZFC}}
\newcommand{\GCH}{\mathsf{GCH}}
\newcommand{\CH}{\mathsf{CH}}
\newcommand{\GA}{\mathsf{GA}}
\newcommand{\GAccc}{{\GA_{\hbox{\scriptsize\sc ccc}}}}
\newcommand{\BA}{\mathsf{BA}}
\newcommand{\HOD}{\mathsf{HOD}}
\newcommand{\restrict}{\upharpoonright}
\newcommand{\df}{\em}
\newcommand{\set}[1]{\{\,{#1}\,\}}
\newcommand{\st}{\mid}
\newcommand{\intersect}{\cap}
\newcommand{\forces}{\Vdash}
\def\<#1>{\langle#1\rangle}
\newcommand{\smalllt}{\mathrel{\mathchoice{\raise2pt\hbox{$\scriptstyle<$}}{\raise1pt\hbox{$\scriptstyle<$}}{\raise0pt\hbox{$\scriptscriptstyle<$}}{\scriptscriptstyle<}}}
\newcommand{\ltdelta}{{{\smalllt}\delta}}
\newtheorem{theorem}{Theorem}

\begin{document}


\address{The College of Staten Island of The City University of New York, Mathematics, 2800 Victory Boulevard, Staten Island, NY 10314, USA}
\address{The Graduate Center of The City University of New York, Mathematics, 365 Fifth Avenue, New York, NY 10016, USA}
\email{jhamkins@gc.cuny.edu, http://jdh.hamkins.org}

\begin{talk}{Joel David Hamkins}{The Ground Axiom}{Hamkins, Joel David}

\noindent Many interesting models of set theory are not
obtainable by nontrivial forcing over an inner model. This
includes, for example, the constructible universe $L$, the
canonical model $L[\mu]$ of a measurable cardinal and many
instances of the core model $K$ (although Schindler has
observed that the least inner model $M_1$ of a Woodin
cardinal actually is a nontrivial forcing extension of an
inner model). To hightlight this phenomenon, my student
Jonas Reitz and I introduced the Ground Axiom, which
asserts that the universe is not a set forcing extension of
any proper inner model.

\newtheorem*{GroundAxiom}{Ground Axiom}
\begin{GroundAxiom}[H, Reitz] The universe is not a forcing extension of any inner model by nontrivial set forcing. Specifically, if $W\ofnoteq V$
is a transitive inner model of $\ZFC$ and $G\of\PP\in W$ is
$W$-generic, then $V\neq W[G]$.
\end{GroundAxiom}

Despite the {\it prima facie} second order nature of this
assertion, the Ground Axiom is actually first order
expressible in the language of set theory.

\begin{theorem}[Reitz, Woodin] The Ground Axiom is first order expressible in the language of $\ZFC$.\label{Theorem.GAisExpressible}
\end{theorem}

This theorem is the starting point of Reitz's dissertation
\cite{Reitz2006:Dissertation}, but an essentially
equivalent assertion was observed independently by Woodin
\cite{Woodin:RecentDevelopmentsOnCH}. Reitz's proof makes
use of ideas arising in Laver's
\cite{Laver:CertainVeryLargeCardinalsNotCreated} recent
result that a ground model is always definable in its
forcing extensions.

\begin{theorem}[Laver] If $V\of V[G]$ is a set forcing extension, then $V$ is a definable class in $V[G]$, using parameters in
$V$.\label{Theorem.GroundIsDefinable}
\end{theorem}

This result was also observed independently by Woodin
\cite{Woodin:RecentDevelopmentsOnCH}. Laver's proof is
connected with my recent theorem showing the extent to
which embeddings in a forcing extension must be lifts of
ground model embeddings.

\newtheorem{KeyDefinition}[theorem]{Key Definition}

\begin{KeyDefinition} \rm\
\begin{enumerate}\item $V\of V[G]$ exhibits {\df $\delta$-covering} if every set of ordinals in $V[G]$ of size less
than $\delta$ is covered by a set of size less than
$\delta$ in $V$. \item $V\of V[G]$ exhibits {\df
$\delta$-approximation} if whenever $A\in V[G]$, $A\of
V$ and $A\intersect a\in V$ for all $a\in V$ with
$|a|^V<\delta$, then $A\in V$.\end{enumerate}
\end{KeyDefinition}

Such forcing extensions are abundant in the large cardinal
literature. Any forcing notion of size less than $\delta$
has $\delta$-approximation and $\delta$-covering. More
generally, any forcing of the form $\PP*\Qdot$, where $\PP$
is nontrivial, $|\PP|<\delta$ and $\forces_\PP\Qdot$ is
$\ltdelta$-strategically closed, exhibits
$\delta$-approximation and $\delta$-covering. Therefore,
such forcing as the Laver preparation or the canonical
forcing of the $\GCH$ exhibit approximation and covering
for many values of $\delta$.

A special case of the main theorem of
\cite{Hamkins:ExtensionsWithApproximationAndCoverProperties}
is:

\begin{theorem}
If\/ $V\of V[G]$ exhibits $\delta$-approximation and
$\delta$-covering, then every ultrapower embedding
$j:V[G]\to M[j(G)]$ above $\delta$ in $V[G]$ is the lift of
an embedding $j\restrict V:V\to M$ definable in
$V$.\label{Theorem.ApproxCoverLift}
\end{theorem}

In particular, $M\of V$ and $j\restrict A\in V$ for all
$A\in V$. The full theorem applies to all sufficiently
closed embeddings, including many types of extender
embeddings. The general conclusion is that extensions with
$\delta$-approximation and $\delta$-covering have no new
large cardinals above $\delta$. The proofs of  Theorems
\ref{Theorem.GroundIsDefinable} and
\ref{Theorem.ApproxCoverLift} make similar and extensive
iterated use of the approximation and cover properties in
their arguments that the respective classes are definable.

Returning to the Ground Axiom, one observes that the
natural models of $\GA$, such as $L$ and $L[\mu]$, exhibit
the $\GCH$ and many other regularity features. Are these a
consequence of the Ground Axiom? The answer is no.

\begin{theorem}[Reitz] If\/ $\ZFC$ is consistent, then $\ZFC+\GA+\neg\CH$ is consistent.\label{Theorem.GA+notCH}
\end{theorem}

The method is flexible and shows that if $\sigma$ is any
$\Sigma_2$ assertion consistent with $\ZFC$, then
$\ZFC+\GA+\sigma$ is consistent. These theorems are proved
by forcing, which is a bit paradoxical as $\GA$ asserts
that the universe is not a forcing extension. Specifically,
resolving the paradox, they are proved by {\it class}
forcing. Using McAloon's
\cite{McAloon1971:ConsistencyResultsAboutOrdinalDefinability}
methods to force strong versions of $V=\HOD$, one codes the
universe into the continuum function, and then $\GA$ holds
with any desired $V_\alpha$ left intact. The hypothesis
$V=\HOD$, however, by itself does not imply $\GA$.
Conversely, at the Set Theory Workshop at the Mathematische
Forshungsinstitut Oberwolfach (0549, December 4-10, 2005),
Woodin suggested a very promising line of argument to show
that the Ground Axiom is consistent with $V\neq\HOD$, which
is now being investigated. Reitz has proved, using large
cardinal indestructibility results, that the Ground Axiom
is consistent with nearly any kind of large cardinal, from
measurable to strong to supercompact and beyond.

\begin{theorem}
If the existence of a supercompact cardinal is consistent
with $\ZFC$, then it is consistent with $\ZFC+\GCH+\GA$.
\end{theorem}

These theorems fit very well into the long-standing set
theoretic program, advanced by Woodin and others, to obtain
the features of the canonical inner models of large
cardinals, but to obtain them by forcing over arbitrary
models of those large cardinals. The Ground Axiom is such a
feature.

Ordinarily, one imagines forcing as a way to reach out into
larger mathematical universes. Here, however, we are
reaching from a given universe down into the possible
ground models of which it is a forcing extension. Given a
model of set theory, perhaps we can strip away a top layer
of forcing and be left with a ground model, a bedrock model
if you will, that is not itself obtainable by forcing from
any smaller inner model. In this case, the original
universe satisfies:

\newtheorem*{BedrockAxiom}{Bedrock Axiom}
\begin{BedrockAxiom}
The universe $V$ is a set forcing extension $V=W[G]$ of an
inner model $W$ of $\ZFC+\GA$.
\end{BedrockAxiom}

The model $W$ is a bedrock model for $V$ in the sense that
it is a minimal ground model for $V$, having no ground
model below it. This axiom is first order expressible for
the same reasons that the Ground Axiom was. Since $V=W$ is
allowed, we have $\GA\implies\BA$. A common feature of the
models of $\GA$ and their forcing extensions, of course, is
that they are all forcing extensions of a model of $\GA$,
and hence themselves models of $\BA$. Are there any other
models? Yes.

\begin{theorem}[Reitz] If\/ $\ZFC$ is consistent, then $\ZFC+\neg\BA$ is consistent. Indeed, if $\sigma$ is any $\Sigma_2$ assertion consistent with
$\ZFC$, then $\ZFC+\BA+\sigma$ is consistent.
\end{theorem}

Perhaps the main open question here is:

\newtheorem{question}[theorem]{Question}
\begin{question}Is the bedrock model unique when it exists?
\end{question}

Several attacks on this question were suggested by various
participants at the Oberwolfach workshop, and a promising
investigation has now ensued.

The theme of current work is to investigate the spectrum of
possible ground models of the universe, the spectrum of
inner models $W$ of which the universe $V$ is a forcing
extension $V=W[G]$. The results above provide a uniform
definition for these ground models $W$ in $V$. By varying
the parameters in this definition, one obtains in effect a
class enumeration of the possible ground models $W$ for
$V$. That is, the class $I$ of parameters $p$ giving rise
to a ground model $W_p$ such that $V$ is a forcing
extension $V=W_p[G_p]$ is definable, and the corresponding
meta-class $\set{W_p\st p\in I}$ of possible ground models
is in effect definable as $\set{\<p,x>\st x\in W_p\And p\in
I}$. Thus, the treatment of the spectrum of possible ground
models is entirely a first order affair of $\ZFC$. Another
theme is to restrict attention to a particular class of
forcing notions, with such axioms as $\GAccc$, which
asserts that the universe is not a nontrivial forcing
extension of an inner model by c.c.c.~forcing. We can
produce models, for example, of $\neg\GA+\GAccc+\sigma$,
for any consistent $\Sigma_2$ assertion $\sigma$; these
models are forcing extensions of an inner model, but are
not obtainable by c.c.c.~forcing. Similar questions and
results abound here.

\end{talk}

\end{document}